\documentclass[10pt]{article}
\usepackage[sectionbib,authoryear]{natbib} 
\usepackage{array}
\usepackage{amssymb}
\usepackage{amsmath}

\usepackage{graphics}
\usepackage{indentfirst}
\usepackage{authblk}
\usepackage{bm}
\usepackage{graphicx}

\newtheorem{thm}{Theorem}

\newtheorem{lem}{Lemma}

\title{Asymptotic power of likelihood ratio tests for high dimensional data}
\author[a,b]{Cheng Wang \thanks{wwcc@mail.ustc.edu.cn}}
\author[b]{Longbing Cao\thanks{LongBing.Cao@uts.edu.au}}
\author[a]{Baiqi Miao\thanks{bqmiao@ustc.edu.cn}}
\affil[a]{Department of Statistics and Finance, University of Science and
Technology of China, Hefei, Anhui 230026, China}
\affil[b]{Advanced Analytics Institute, University of Technology Sydney, New
South Wales 2007, Australia}
\date{}

\begin{document}
\maketitle

\begin{abstract}
This paper considers the asymptotic power of likelihood ratio test (LRT) for the identity test when the dimension $p$ is large compared to the sample size $n$. The asymptotic distribution of LRT under alternatives is given and an explicit expression of the power is derived when $\rm{tr}(\Sigma_p)-\log|\Sigma_p|-p$ tends to a constant. A simulation study is carried out to compare LRT with other tests.  All these studies show that LRT is a powerful test to detect eigenvalues around zero. 

\end{abstract}

% \linenumbers

\section{Introduction}
In multivariate analysis for high dimensional data, testing the structure of population covariance matrices is an important problem. See, for example, \cite{johnstone2001distribution}, \cite{lw}, \cite{srivastava2005some},  \cite{schott2006high},  \cite{Chen10}, \cite{cai2011limiting} and \cite{li2012two}, among many others. Let $X_1,\cdots,X_n$ be $n$ independent and identically distributed (i.i.d.) $p$-variate random vectors following a multivariate normal distribution $N_p(\mu,\Sigma_p)$ where $\mu$ is the mean vector and $\Sigma_p$ is the covariance matrix. In many studies, a hypothesis
test of significant interest is to test
\begin{eqnarray} \label{h0}
H_0:~\Sigma_p=I_p~~vs~~H_1:~\Sigma_p \neq I_p,
\end{eqnarray}
where $I_p$ is the $p$-dimensional identity matrix. Note that the identity matrix  in (\ref{h0}) can be replaced by any other
positive definite matrix $\Sigma_0$ through multiplying the data by $\Sigma_0^{-1/2}$.

A natural approach to test (\ref{h0}) is to conduct estimations for some distance measures between 
$\Sigma_p$ and $I_p$ and there are two types of measures which are widely used in literature. The first is based on the likelihood function:
\begin{eqnarray}
L_l(\Sigma_p)=\rm{tr}(\Sigma_p)-\log{|\Sigma_p|}-p,
\end{eqnarray} 
and the second is based on quadratic loss function:
\begin{eqnarray}
L_q(\Sigma_p)=\rm{tr}(\Sigma_p-I_p)^2.
\end{eqnarray}
Each of these is 0 when $\Sigma_p=I_p$ and positive when $\Sigma_p \neq I_p$.
For $L_l(\Sigma_p)$, it is referred to the likelihood ratio test (LRT) and the classic LRT for fixed $p$ and
large $n$ can be found in \cite{And}. For high dimensional data ($p$ is large), the failure of classical LRT was firstly observed by \cite{dempster1958high} and later in a pioneer work by \cite{Bai09}, authors
proposed corrections to LRT when $p/n \to c \in (0,1)$ and $\mu=0$. Successive
works included \cite{jiang2012likelihood} which extended the results of \cite{Bai09} to Gaussian data
with general $\mu$ and our work \cite{wang2012identity} where we studied the LRT for general $\mu$ and 
non-Gaussian data.  For the quadratic loss function $L_q(\Sigma_p)$, there are many works since the
seminal paper \cite{nagao1973some}. These included \cite{lw}, \cite{Birke05}, \cite{srivastava2005some}, \cite{Chen10}
and \cite{cai2012optimal}. Other works which considered question (\ref{h0}) are referred to \cite{johnstone2001distribution}, \cite{cai2011limiting} and \cite{onatski2011asymptotic}. 

The existing results about LRT \citep{Bai09,jiang2012likelihood, wang2012identity} have only derived asymptotic null distribution and we know little about the asymptotic point-wise power of LRT under the alternative hypothesis. Recently, \cite{onatski2011asymptotic} studied the asymptotic power of several tests including LRT under the special alternative of rank one perturbation to the identity matrix as both $p$ and $n$ go to infinity. \cite{cai2012optimal}
investigated the testing problem (\ref{h0}) in the high-dimensional settings from a minimax point of view. The results in \cite{onatski2011asymptotic} and \cite{cai2012optimal} showed that LRT was a sub-optimal test when there was a rank one perturbation to the identity matrix . 

In this work, we will consider the power of LRT under general alternatives. The asymptotic distribution of LRT will be studied when $\Sigma_p \neq I_p$ and an explicit expression of the power will also be derived when $L_l(\Sigma_p)$ tends to a constant. From these results, we find that in relation to LRT it is not fair that \cite{onatski2011asymptotic} and \cite{cai2012optimal} only focused on the alternatives whose eigenvalues were larger than 1. Furthermore, our results show that LRT is powerful to detect eigenvalues around zero. Simulations will also be conducted to compare LRT with two tests based on $L_q(\Sigma_p)$ \citep{Chen10, cai2012optimal}.       
   
The paper is structured as follows: Section 2 introduces the basic data
structure and establishes the asymptotic power of LRT while Section 3 reports simulation studies. All the proofs are presented in the Appendix.
\section{Main Results}
To relax the Gaussian assumptions, we assume that the observations $X_1,\cdots,X_n$ satisfy a multivariate model
\citep{Chen10}
\begin{eqnarray} \label{data}
X_i=\Sigma_p^{1/2} Y_i+\mu,~for~i=1,\cdots,n
\end{eqnarray}
where $\mu$ is a $p$-dimensional constant vector and the entries of
$\mathcal{Y}_n=(Y_{ij})_{p \times n}=(Y_1,\cdots,Y_n)$ are i.i.d. with $E
Y_{ij}=0$, $EY^2_{ij}=1$ and $EY^4_{ij}=3+\Delta$. The sample covariance matrix is defined as
\begin{eqnarray*}
S_n=\frac{1}{n-1} \sum_{k=1}^n (X_k-\bar{X}) (X_k-\bar{X})',
\end{eqnarray*}
where $\bar{X}=\frac{1}{n} \sum_{k=1}^n X_k$.

Writing $y_n=p/n<1$, the LRT statistic is defined as
\begin{eqnarray} \label{lrt}
L_n=\frac{1}{p} \rm{tr}(S_n)-\frac{1}{p} \log{|S_n|}-1-d(y_n)
\end{eqnarray}
where $\rm{tr}$ denotes the trace and $d(x)=1+(1/x-1) \log{(1-x)},~0<x<1$. Under the null hypothesis, \cite{wang2012identity} derived the following asymptotic normality of $L_n$ by using random matrix theories.
\begin{thm}[Theorem 2.1 of \cite{wang2012identity}] \label{thm1}
When $\Sigma_p=I_p$ and $y_n=p/y \to y \in(0,1)$,
\begin{eqnarray*}
\frac{p L_n-\mu_n}{\sigma_n} \stackrel{D}{\rightarrow} N(0,1),
\end{eqnarray*}
where $\mu_n=y_n (\Delta/2-1) -3/2 \log(1-y_n)$, $\sigma_n^2=-2 y_n- 2\log{(1-y_n)}$ and $\stackrel{D}{\rightarrow}$ denotes convergence in distribution.
\end{thm}

When $X_1, \cdots, X_n$ be i.i.d. multivariate normal distributions $N_p(\mu,\Sigma_p)$ where $\Delta=0$, \cite{jiang2012likelihood} derived a similar result as Theorem \ref{thm1} by using the Selberg integral and they also considered the special situation where $p/n \to 1$. Based on the asymptotic normality under the respective null
hypothesis, an asymptotic level $\alpha$ test based on $L_n$ is given by
\begin{eqnarray} \label{test}
\phi=I(\frac{p L_n-\mu_n}{\sigma_n}>z_{1-\alpha}),
\end{eqnarray}
where $I(\cdot )$ is the indicator function, and $z_{1-\alpha}$ denotes the $100 \times (1-\alpha)$th percentile of the
standard normal distribution. In the following theorem, we establish the convergence of $L_n$ under the alternative $\Sigma_p \neq I_p$.
\begin{thm} \label{thm2}
When $\rm{tr}(\Sigma_p-I_p)^2/p \to 0$ and $y_n=p/y \to y \in(0,1)$,
\begin{eqnarray*}
\frac{p L_n-L_l(\Sigma_p)-\mu_n}{\sigma_n} \stackrel{D}{\rightarrow} N(0,1),
\end{eqnarray*}
where $\mu_n=y_n (\Delta/2-1) -3\log(1-y_n)/2$ and $\sigma_n^2=-2 y_n- 2\log{(1-y_n)}$.
\end{thm}
 
 In particular, when $L_l(\Sigma_p)$ tends to a constant, we have the following result.
\begin{thm} \label{thm3}
When $ L_l(\Sigma_p) \to b \in (0,\infty)$ and $y_n=p/y \to y \in(0,1)$,
\begin{eqnarray*}
\lim_{n \to \infty} P_{\Sigma_p}(\phi~rejects~H_0)=1-\Phi(z_{1-\alpha}-\frac{b}{\sqrt{-2 y-2 \log{(1-y)}}}),
\end{eqnarray*}
where $\Phi(\cdot)$ is the cumulative distribution function of the standard normal distribution.
\end{thm}
It can be seen from Theorems \ref{thm2} and \ref{thm3} that the expression
\begin{eqnarray}
1-\Phi(z_{1-\alpha}-\frac{L_l(\Sigma_p)}{\sigma_n}),
\end{eqnarray}
gives good approximation to the power of the test in (\ref{test}) until the power is extremely
close to 1. In particular, when $L_l(\Sigma_p)$ is large, the power of the test $\phi$ will be close to 1 and it is hard for $\phi$ to distinguish between the two hypotheses if $L_l(\Sigma_p)$ tends to zero. 

To derive the asymptotic power, a special covariance matrix was used in \cite{cai2012optimal} and \cite{onatski2011asymptotic} as follows
\begin{eqnarray} \label{alsigma}
\Sigma_p^{*}=I_p+h\sqrt{\frac{p}{n}} v v',
\end{eqnarray} 
where $h$ is a constant and $v$ is an arbitrarily fixed unit vector. For this special spiked matrix  \citep{johnstone2001distribution}, the true $\Sigma_p$ has a perturbation in a single unknown direction and in \cite{cai2012optimal} and \cite{onatski2011asymptotic}, the authors focused on situations where $h>0$ that is 
one eigenvalue of $\Sigma_p$ is larger than 1 while others are still unitary. Here, by Theorem \ref{thm3}, we know that the asymptotic power of the LRT test (\ref{lrt}) is  
\begin{eqnarray*}
1-\Phi(z_{1-\alpha}-\frac{h \sqrt{y}-\log(1+h \sqrt{y})}{\sqrt{-2 y-2 \log{(1-y)}}}).
\end{eqnarray*}
Therefore, compared with the tests based on $L_q(\Sigma_p)$ \citep{lw, Chen10, cai2012optimal} whose power for $\Sigma_p^*$ is $1-\Phi(z_{1-\alpha}-h^2/2)$, LRT is more sensitive to small eigenvalues ($h<0$), not any bigger than one ($h>0$). In particular, when $1+h \sqrt{y}$ is close to 0 that is $\Sigma_p^*$ has a very small eigenvalue, the power will tend to 1. A numerical experiment will be conducted in the next section to show the performances of LRT (\ref{test}) and the tests based on $L_q(\Sigma_p)$.    
\section{Simulations}
In this section, we conduct several simulation studies to compare the power of the LRT in (\ref{test}) with that of the 
tests based on $L_q(\Sigma_p)$. When $X_1,\cdots,X_n$ i.i.d from $N_p(0,\Sigma_p)$, \cite{cai2012optimal} proposed an estimator of $L_q(\Sigma_p)$ as
\begin{eqnarray}
T_{1,n}=\frac{2}{n(n-1)} \sum_{1\leq i<\leq n} h(X_i,X_j),
\end{eqnarray}
where $h(X_1,X_2)=(X_1'X_2)^2-(X_1'X_1+X_2'X_2)+p$ and an asymptotic level $\alpha$ test based on $T_{1,n}$ is given by
\begin{eqnarray}
\phi_1=I(T_{1,n}>z_{1-\alpha}2 \sqrt{\frac{p(p+1)}{n(n-1)}}).
\end{eqnarray}
Similarly, for general data structure, \cite{Chen10} gave an estimator for $L_q(\Sigma_p)$ as
\begin{eqnarray*}
T_{2,n}=&&\frac{1}{P_n^2}\sum_{i,j}^*(X_i'X_j)^2-\frac{2}{P_n^3} \sum_{i,j,k}^*X_i'X_jX_j'X_k\\
&&+\frac{1}{P_n^4}\sum_{i,j,k,l}^*X_i'X_jX_k'X_l-2 tr(S_n)+p,
\end{eqnarray*}
where $P_n^r=n!/(n-r)!$ and $\sum^*$ denotes summation over
mutually different indices. And the level $\alpha$ test is  
\begin{eqnarray}
\phi_2=I(\frac{n}{2 p}T_{2,n}>z_{1-\alpha}).
\end{eqnarray}
To evaluate the power of the tests, the alternative population covariance matrix will be set as $\Sigma^*=diag(\rho,1,\cdots,1)$ where $\rho$ will range from $0.01$ to 4. Here, we only focus on this simple matrix to investigate the features of each test and simulations for more general alternatives can be found in our previous work \cite{wang2012identity}. All the results are based on $10^4$ replications.

\begin{figure}
\centering
\includegraphics[scale=0.45]{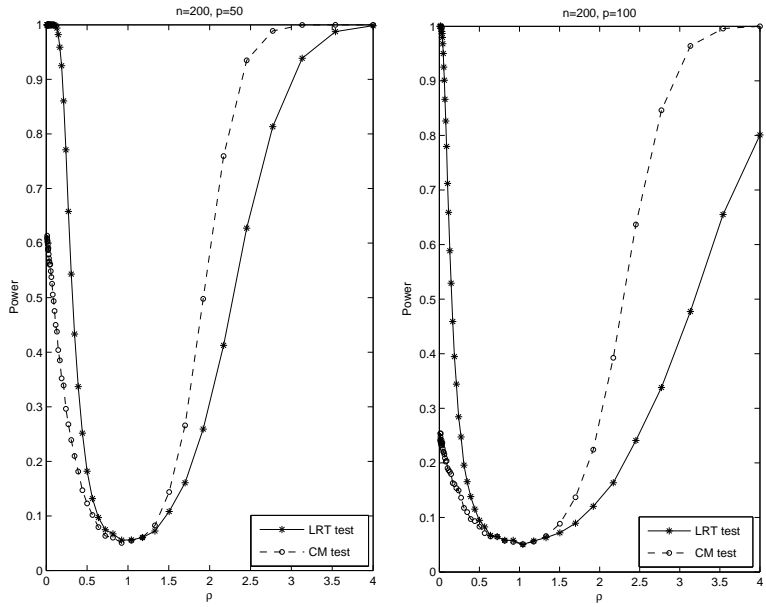}
\caption{Empirical performances of LRT and CM tests for Gaussian data with known means.}
\label{fig2}
\end{figure}

Figure \ref{fig2} reports empirical powers of the LRT and the test of \cite{cai2012optimal} (CM test) for $N_p(0,\Sigma_p)$ distributions. We observe from Figure \ref{fig2} that LRT had a better performance for $\rho<1$ while the CM test performed better
if $\rho>1$. When $\rho$ is around 1 that is the true $\Sigma_p$ is very close to an identity matrix, both tests had similar empirical sizes which were quite close to the nominal 5\%. For the fixed sample size $n=200$, when $p$ was increased from 
$50$ to $100$, the performances of both tests became poor which could be understood as the estimators get worse as $p$ is increased for the fixed sample size.  

\begin{figure}
\centering
\includegraphics[scale=0.45]{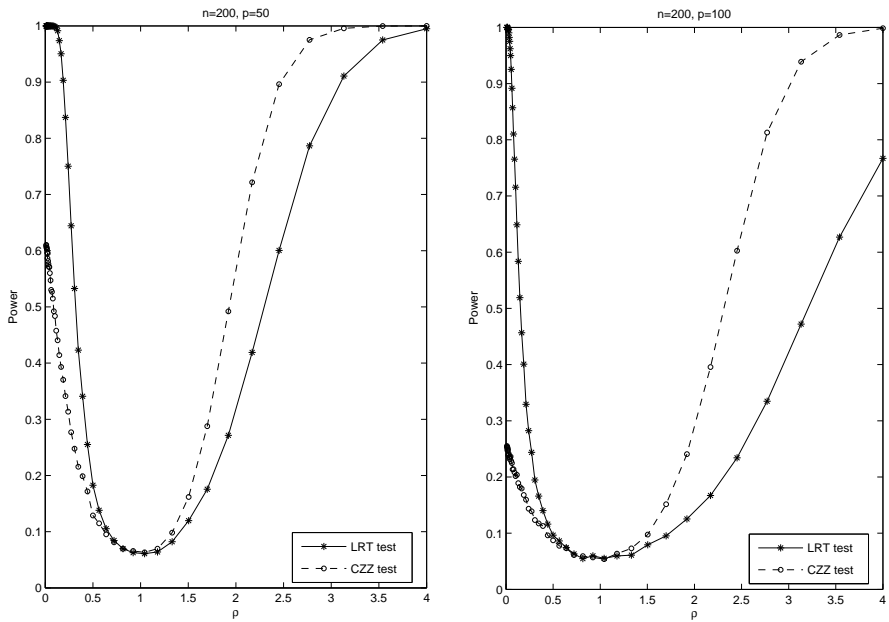}
\caption{Empirical performances of LRT and CZZ tests for non-Gaussian data with unknown means}
\label{fig3}
\end{figure}

For general data, due to the CM test is only applicable for $N_p(0,\Sigma_p)$, we conduct comparisons between LRT and the test of \cite{Chen10} (CZZ test) where $Y_{ij}$ comes from Gamma(4, 0.5) distribution in data (\ref{data}). Figure \ref{fig3} reports the empirical powers of the LRT and CZZ tests and the conclusions follow very
similar patterns to those of Figure \ref{fig2} which shows that LRT is a powerful test to detect eigenvalues around zero.

\section{Appendix}
\subsection{Proof of Theorem \ref{thm2}}
Writing 
\begin{eqnarray}
B_n=\frac{1}{n-1} \sum_{k=1}^n (Y_k-\bar{Y})(Y_k-\bar{Y})',
\end{eqnarray}
where $\bar{Y}=\sum_{k=1}^{n} Y_k$. Noting $S_n=\Sigma_p^{1/2} B_n \Sigma_p^{1/2}$, we have
\begin{eqnarray*}
L_n&=&\frac{1}{p} \rm{tr}(S_n)-\frac{1}{p} \log{|S_n|}-1-d(y_n)\\
&=& \frac{1}{p} \rm{tr}(\Sigma_p B_n)-\frac{1}{p} \rm{tr}(B_n)-\frac{1}{p} \log{|\Sigma_p|}+BL_{n}\\
&=& BL_{n}+\frac{1}{p} L_l(\Sigma_p)+ \frac{1}{p} \rm{tr}((\Sigma_p-I_p) B_n)-\frac{1}{p} \rm{tr}(\Sigma_p-I_p),
\end{eqnarray*}
where $BL_n=\frac{1}{p} \rm{tr}(B_n)-\frac{1}{p} \log{|B_n|}-1-d(y_n)$. By Theorem \ref{thm1},
\begin{eqnarray*}
\frac{p BL_n-\mu_n}{\sigma_n} \stackrel{D}{\rightarrow} N(0,1).
\end{eqnarray*}
Therefore, to prove Theorem \ref{thm2}, it is enough to show
\begin{eqnarray*}
\epsilon_n:=\rm{tr}((\Sigma_p-I_p) B_n)- \rm{tr}(\Sigma_p-I_p)=o_p(1).
\end{eqnarray*}
Rewriting $B_n$ as 
\begin{eqnarray}
B_n=\frac{1}{n} \sum_{k=1}^n Y_k Y_k'-\frac{1}{n(n-1)}\sum_{i \neq j}Y_iY_j',
\end{eqnarray}
we can get $E[\epsilon_n]=0$ and by Proposition A.1 of \cite{Chen10},
\begin{eqnarray*}
E[\epsilon_n^2]&=&E[(\frac{1}{n} \sum_{k=1}^n Y_k'(\Sigma_p-I_p) Y_k-\frac{1}{n(n-1)}\sum_{i \neq j}Y_i'(\Sigma_p-I_p)Y_j)^2]-(\rm{tr}(\Sigma_p-I_p))^2\\
&=& \frac{2}{n-1} \rm{tr}((\Sigma_p-I_p)^2)+ \frac{\Delta}{n} tr((\Sigma_p-I_p)\circ (\Sigma_p-I_p))\\
&\leq & \frac{2+\Delta}{n-1} \rm{tr}((\Sigma_p-I_p)^2)
\end{eqnarray*}
where $\circ$ denotes Hadamard product. Above all, when $p/n  \to y$ and $\rm{tr}((\Sigma_p-I_p)^2)/p \to 0$, we come to  $\epsilon_n=o_p(1)$.

The proof is completed. 
\subsection{Proof of Theorem \ref{thm3}}
To prove the theorem, we need some inequalities.
\begin{lem} \label{lem1}
For any $x>0$, 
\begin{itemize}
\item[(1)] If $0<x \leq 1$, $(x-1)^2 \leq 2 (x-1-\log{x})$;
\item[(2)] If $1<x<M$, $(x-1)^2 \leq 2 M (x-1-\log{x})$. 
\end{itemize} 
\end{lem}

Since $ L_l(\Sigma_p) \to b \in (0,\infty)$, for large enough $p$, we will always have $D_1(\Sigma_p)<2 b$.
Denoting the eigenvalues of $\Sigma_p$ as $\lambda_{1,p}, \cdots, \lambda_{p,p}$, 
\begin{eqnarray*}
L_l(\Sigma_p)=\sum_{k=1}^p (\lambda_{k,p}-\log{\lambda_{k,p}}-1)<2b,
\end{eqnarray*}
which implies that there is a constant $c_0=c_0(b)>1$ satisfying 
\begin{eqnarray*}
c_0^{-1} \leq \lambda_{1,p}, \cdots, \lambda_{p,p} \leq c_0.
\end{eqnarray*} 
By Lemma \ref{lem1}, $(\lambda_{k,p}-1)^2 \leq 2 c_0(\lambda_{k,p}-\log{\lambda_{k,p}}-1)$ which means
\begin{eqnarray*}
L_q(\Sigma_p)=\sum_{k=1}^p (\lambda_{k,p}-1)^2\leq 2c_0 L_l(\Sigma_p) < 4 c_0 b=o(p).
\end{eqnarray*}
By Theorem \ref{thm2} and Slutsky's lemma, 
\begin{eqnarray*}
\frac{p L_n-\mu_n}{\sigma_n}-\frac{b}{\sqrt{-2 y-2 \log{(1-y)}}}\stackrel{D}{\rightarrow} N(0,1).
\end{eqnarray*}
Now, we can calculate the power of the test
\begin{eqnarray*}
& &P_{\Sigma_p}(\phi~rejects~H_0)=P(\frac{p L_n-\mu_n}{\sigma_n}>z_{1-\alpha})\\
&=&P(\frac{p L_n-\mu_n}{\sigma_n}-\frac{b}{\sqrt{-2 y-2 \log{(1-y)}}}>z_{1-\alpha}-\frac{b}{\sqrt{-2 y-2 \log{(1-y)}}})\\
&\to & 1-\Phi(z_{1-\alpha}-\frac{b}{\sqrt{-2 y-2 \log{(1-y)}}}).
\end{eqnarray*}
The proof is completed. 

\section*{Acknowledgement}

We thank Dr. Guangming Pan and Dr. Zongming Ma for their helpful discussions and suggestions. The research of Cheng Wang and Baiqi Miao was partly supported by NSF of China Grands No. 11101397 and 71001095. Longbing Cao's research was supported in part by the Australian Research Council Discovery Grant DP1096218 and the Australian
Research Council Linkage Grant LP100200774.

\bibliographystyle{natbib}
\bibliography{cite}

\begin{thebibliography}{}

\bibitem[Anderson(2003)Anderson]{And}
Anderson, T. (2003).
\newblock {\em An introduction to multivariate statistical analysis.}
\newblock Hoboken, NJ:Wiley.

\bibitem[Bai {\em et~al.}(2009)Bai, Jiang, Yao, and Zheng]{Bai09}
Bai, Z., Jiang, D., Yao, J., and Zheng, S. (2009).
\newblock Corrections to {LRT} on large-dimensional covariance matrix by rmt.
\newblock {\em The Annals of Statistics\/}, {\bf 37}(6B), 3822--3840.

\bibitem[Birke and Dette(2005)Birke and Dette]{Birke05}
Birke, M. and Dette, H. (2005).
\newblock A note on testing the covariance matrix for large dimension.
\newblock {\em Statistics \& Probability Letters\/}, {\bf 74}(3), 281--289.

\bibitem[Cai and Jiang(2011)Cai and Jiang]{cai2011limiting}
Cai, T. and Jiang, T. (2011).
\newblock Limiting laws of coherence of random matrices with applications to
  testing covariance structure and construction of compressed sensing matrices.
\newblock {\em The Annals of Statistics\/}, {\bf 39}(3), 1496--1525.

\bibitem[Cai and Ma(2012)Cai and Ma]{cai2012optimal}
Cai, T. and Ma, Z. (2012).
\newblock Optimal hypothesis testing for high dimensional covariance matrices.
\newblock {\em arXiv:1205.4219\/}.

\bibitem[Chen {\em et~al.}(2010)Chen, Zhang, and Zhong]{Chen10}
Chen, S., Zhang, L., and Zhong, P. (2010).
\newblock Tests for high-dimensional covariance matrices.
\newblock {\em Journal of the American Statistical Association\/}, {\bf
  105}(490), 810--819.

\bibitem[Dempster(1958)Dempster]{dempster1958high}
Dempster, A. (1958).
\newblock A high dimensional two sample significance test.
\newblock {\em The Annals of Mathematical Statistics\/}, {\bf 29}(4),
  995--1010.

\bibitem[Jiang {\em et~al.}(2012)Jiang, Jiang, and Yang]{jiang2012likelihood}
Jiang, D., Jiang, T., and Yang, F. (2012).
\newblock Likelihood ratio tests for covariance matrices of high-dimensional
  normal distributions.
\newblock {\em Journal of Statistical Planning and Inference\/}.

\bibitem[Johnstone(2001)Johnstone]{johnstone2001distribution}
Johnstone, I. (2001).
\newblock On the distribution of the largest eigenvalue in principal components
  analysis.
\newblock {\em The Annals of Statistics\/}, {\bf 29}(2), 295--327.

\bibitem[Ledoit and Wolf(2002)Ledoit and Wolf]{lw}
Ledoit, O. and Wolf, M. (2002).
\newblock Some hypothesis tests for the covariance matrix when the dimension is
  large compared to the sample size.
\newblock {\em The Annals of Statistics\/}, pages 1081--1102.

\bibitem[Li and Chen(2012)Li and Chen]{li2012two}
Li, J. and Chen, S. (2012).
\newblock Two sample tests for high-dimensional covariance matrices.
\newblock {\em The Annals of Statistics\/}, {\bf 40}(2), 908--940.

\bibitem[Nagao(1973)Nagao]{nagao1973some}
Nagao, H. (1973).
\newblock On some test criteria for covariance matrix.
\newblock {\em The Annals of Statistics\/}, pages 700--709.

\bibitem[Onatski {\em et~al.}(2011)Onatski, Moreira, and
  Hallin]{onatski2011asymptotic}
Onatski, A., Moreira, M., and Hallin, M. (2011).
\newblock Asymptotic power of sphericity tests for high-dimensional data.
\newblock {\em Manuscript\/}.

\bibitem[Schott(2006)Schott]{schott2006high}
Schott, J. (2006).
\newblock A high-dimensional test for the equality of the smallest eigenvalues
  of a covariance matrix.
\newblock {\em Journal of Multivariate Analysis\/}, {\bf 97}(4), 827--843.

\bibitem[Srivastava(2005)Srivastava]{srivastava2005some}
Srivastava, M. (2005).
\newblock Some tests concerning the covariance matrix in high dimensional data.
\newblock {\em J. Japan Statist. Soc\/}, {\bf 35}(2), 251--272.

\bibitem[Wang {\em et~al.}(2012)Wang, Yang, Miao, and Cao]{wang2012identity}
Wang, C., Yang, J., Miao, B., and Cao, L. (2012).
\newblock On identity tests for high dimensional data using rmt.
\newblock {\em arXiv:1203.3278\/}.

\end{thebibliography}
\end{document}